\documentclass[11pt]{article}

\usepackage{amsmath}
\usepackage{amsfonts}
\usepackage{mathrsfs}
\usepackage[latin1]{inputenc}
\usepackage{color}

 \newtheorem{thm}{Theorem}[section]

 \newtheorem{prop}[thm]{Proposition}

 \newtheorem{rem}[thm]{Remark}

\newcommand{\eps}{\varepsilon}
\newcommand{\R}{\mathbb R}
\newcommand{\K}{\mathbb K}
\newcommand{\N}{\mathbb N}
\newcommand{\lraup}{\relbar\joinrel\relbar\joinrel\relbar\joinrel\rightharpoonup}
\newcommand{\cchi}{{\mbox{\raisebox{3pt}{$\chi$}}}}

\title{\huge\bf Quasivariational solutions for first order quasilinear equations with gradient constraint}

\author{Jos\'e Francisco Rodrigues$\quad\qquad$ Lisa
Santos}

\date{}

\begin{document}

\maketitle

\begin{abstract}
We prove the existence of solutions for an evolution quasi-variational inequality with a first order quasilinear operator and a variable convex set,
 which is characterized by a constraint on the absolute value of the gradient that depends on the solution itself. The only required assumption on the
  nonlinearity of this constraint is its continuity and positivity. The
method relies on an appropriate parabolic regularization and suitable {\em a priori} estimates. We obtain also the existence of stationary solutions,
by studying the asymptotic behaviour in time. In the variational case, corresponding to a constraint independent of the solution, we also give
uniqueness results.
\end{abstract}

\baselineskip 16pt

\vspace{5mm}


\section{Introduction}

The initial-boundary value problem (IBVP) for scalar equations of first order
\begin{equation}\label{ibvp}
\partial_t u-\nabla\cdot\boldsymbol\Phi(u)=f(u)
\end{equation}
in bounded sets $\Omega\subset\R^N$ with a smooth boundary $\partial\Omega$ and time $t>0$ may, in general, be not well posed, even for smooth
flux function $\boldsymbol\Phi=\boldsymbol\Phi(x,t,u)$, source term $f=f(x,t,u)$ and initial-boundary smooth data (here we denote by $\partial_t$ the
partial derivative in order to $t$ and by $\nabla\cdot\,$ the divergence in $x\in\Omega$). In the pioneer paper \cite{Bardos_leRouxNedelec1979}, Bardos, Leroux
and N\'ed\'elec, using the method of vanishing viscosity, extended to bounded domains the notion of entropy solution, obtaining their existence and
uniqueness in the BV framework. Dealing with  data merely in $L^\infty$, Otto (\cite{Otto1996}, see also \cite{MalekNecasRokytaRuzicka1996} or
 \cite{Dafermos2010}) has shown the well posedeness  of the IBVP introducing an appropriate weak formulation. Recent results on $L^\infty$ entropy
 solutions were obtained in \cite{Martin2007} and in \cite{CocliteKarlsenKwon2009}, where the delicate boundary trace question is analyzed and more
 references can be found.
 The need of weak formulations  is justified, not only by the possibility of shock fronts but also, in general, by the
 boundary layers introduced by the vanishing viscosity and the impossibility of prescribing Dirichlet data for solutions of (\ref{ibvp}) on the
 whole boundary $\partial\Omega$.

 Here we are interested in solutions of (\ref{ibvp}) with the spatial gradient globally bounded, which not only prevents the existence of shocks but
 also allows the prescription of the data on the whole $\partial\Omega$. This will be done by imposing a gradient constraint to the IBVP and has a motivation
 from critical-state problems.

It is known that several critical-state problems arising in different physical models, such as in elastic-plastic deformations, sandpile growth,
magnetization of type-II superconductors or formation of networks of lakes and rivers, are formulated with convex constraints on first order
derivatives (see, for instance, \cite{Bean1964}, \cite{DuvautLions1972} or \cite{Prigozhin1994} and their references). More than one
approach has been considered leading, in particular, to several works with gradient constraint for the scalar function $u=u(x,t)$,
\begin{equation}\label{guu}
|\nabla u|\le G(u),
\end{equation}
where $G=G(x,u)>0$ is a given threshold function. For the  elastic-plastic torsion problem, the variational inequality formulation of the elliptic
 problem  for the case $G\equiv 1$
is well-known (see \cite{Ting1967},  \cite{BrezisStamppachia1968} or  \cite{Rodrigues1987}, for instance) and has been
extended to elliptic quasivariational inequalities in \cite{KunzeRodrigues2000}, when $G=G(u)$.

The gradient constraint (\ref{guu}), for each $t$, splits the domain $\Omega$ into two regions
\begin{equation*}
\Lambda(t)=\big\{x\in\Omega:\big|\nabla u(x,t)\big|<G(x,u(x,t))\big\},
\end{equation*}
\begin{equation*}
I(t)=\big\{x\in\Omega:\big|\nabla u(x,t)\big|=G(x,u(x,t))\big\},
\end{equation*}
which are, formally, separated by an  unknown free boundary $\partial I(t)\cap\Omega=\partial\Lambda(t)\cap\Omega$. In this work we are
interested in imposing the constraint (\ref{guu}) to solutions of first order quasilinear equations of the type (\ref{ibvp}). In fact, we look for
a function $u=u(x,t)$ that solves the equation (\ref{ibvp}) in the {\em a priori} unknown region $\Lambda(t)$, i.e., if the gradient threshold
is not attained, while, in the complementary region $I(t)$, $u$ solves the Hamilton-Jacobi equation $|\nabla u|=G(u)$, with the prescribed
boundary and initial values $u=0$ on $\partial\Omega\times(0,T)$ and $u(x,0)=u_0(x)$ on $\Omega$, respectively.

Similarly to the elastic-plastic torsion problem, as it was shown by Br\'ezis \cite{Brezis1972} (see also \cite{Rodrigues1987}, p.266), we may
introduce a formal Lagrange multiplier $\lambda=\lambda(x,t)$ associated with the constraint (\ref{guu}) and such that
\begin{equation}\label{lm}
\lambda\ge 0,\qquad\lambda(G(u)-|\nabla u|)=0,
\end{equation}
and
\begin{equation}\label{eq}
\partial_t u-\nabla\cdot\boldsymbol\Phi(u)-\nabla\cdot(\lambda\nabla u)=f(u),
\end{equation}
where now these conditions hold in the whole space-time cylinder $Q=\Omega\times(0,T)$ in some sense. Then, introducing
for each $t\in(0,T)$ the convex set
\begin{equation*}
\K_{G(u(t))}=\big\{v\in H^1_0(\Omega):\big|\nabla v\big|\le G(u(t))\text{ in }\Omega\big\}
\end{equation*}
and integrating by parts (\ref{eq}) with $v-u(t)$, for $v\in\K_{G(u(t))}$ we observe that, for $u=u(t)$, (\ref{lm}) yields
\begin{equation*}\int_\Omega\lambda\nabla u\cdot\nabla(v-u)\le\int_\Omega\lambda\big(\big|\nabla u\big|\,\big|\nabla v\big|
-\big|\nabla u\big|^2\big)\le\int_\Omega\lambda\big(G(u)-\big|\nabla u\big|\big)\big|\nabla u\big|=0
\end{equation*}
and then $u(t)\in\K_{G(u(t))}$ satisfies the quasivariational inequality, for a.e. $t\in(0,T)$,
\begin{equation}\label{6}
\int_\Omega\partial_tu(v-u)+\int_\Omega\boldsymbol\Phi(u)\cdot\nabla(v-u)\ge\int_\Omega f(u)(v-u),\qquad \forall v\in\K_{G(u(t))},
\end{equation}
with the prescribed initial condition $u(0)=u_0$.

This is a natural formulation for general first order quasilinear scalar operators under the constraint (\ref{guu}) and
it leads to a new class of problems that were not considered in the classical references of quasivariational inequalities
\cite{BaiocchiCapelo1984} or \cite{BensoussanLions1982}. However, it contains as a special case the problem considered in \cite{BarrettPrigozhin2010}
motivated by a critical-state model in superconductivity. It was also considered in the setting of parabolic operators in
\cite{Santos1991,Santos2002},
in the variational inequality case, i.e., when $G$ does not depend on the solution $u$, and in a quasivariational case by the authors
in \cite{RodriguesSantos2000}. We note that Prigozhin has considered a similar evolutionary variational inequality in superconductivity
\cite{Prigozhin1996} and for a sandpile model in \cite{Prigozhin1996-3}, corresponding essentially to the particular case $\boldsymbol\Phi\equiv 0$,
 $f\equiv$\,constant$>0$
and $G\equiv 1$. This last model, which can also be interpreted as an evolution governed by the multivalued differential associated with the constraint
$|\nabla v|\le 1$, has also been considered by other authors (see, for instance, \cite{AronsonEvansWu1996} or \cite{CannarsaCardaliaguetSinestrari2009} and their references).
Recent results for parabolic problems with gradient constraint have been obtained in \cite{AzevedoSantos2010} and \cite{Kenmochi2012}.

Under the natural assumptions on the smoothness of $\boldsymbol\Phi$ and $f$ and only positivity and continuity on $u$ of $G$, we prove the existence of
a globally continuous solution $u$ to the quasivariational inequality (\ref{6}), with a bounded measure time derivative $\partial_t u$.
For the special case when $G=G(x)>0$ is merely bounded we show the uniqueness of solution $\vartheta(t)$ of the variational inequality (\ref{6}) in $\K_G$
and the $L^2$ integrability of $\partial_t\vartheta$. Under appropriate assumptions on the data we also prove the asymptotic behaviour of the
evolutionary solutions as $t\rightarrow\infty$. We obtain, in particular, the existence of a solution $u^\infty$ to the stationary quasivariational
inequality (\ref{6}), when $\boldsymbol\Phi$ is time independent, $f$ is strictly decreasing in $u$ and $G$ is bounded in $u$ or $\nabla\cdot
\boldsymbol\Phi+f$ is compatible with the weak maximum principle. In the variational inequality case, if $f^\infty$ is  strictly decreasing
in $u$, we have uniqueness of the stationary solution $\vartheta^\infty$ and we obtain an estimate of the asymptotic stabilization of $\vartheta(t)
\longrightarrow\vartheta^\infty$ as $t\rightarrow\infty$, first in $L^1$ and then also in Hölder spaces.

The precise assumptions and main results are stated in Section 2. We use the method of vanishing viscosity and constraint penalization in two
consecutive steps, by considering an approximating quasilinear parabolic problem with convection, extending the approach of \cite{RodriguesSantos2000}.
The essential {\em a priori} estimates are established in Section 3 and the passage to the limits in Section 4. In the intermediate step,
 we show the existence of a solution to the new parabolic quasivariational inequality (\ref{iqvdelta}), generalizing a result of
 \cite{RodriguesSantos2000} in the case $p=2$.

 The existence of unique solution of the variational inequality case is proven in Section 5, extending the techniques used in \cite{Santos1991} and
 \cite{Santos2002} for the parabolic case. We also give an indirect proof of the existence  of the quasivariational stationary solution in the final Section 6. There
  we prove the asymptotic behaviour as $t\rightarrow\infty$, extending the estimate of \cite{RodriguesSantos2000} in $\partial_t u$
 for the quasivariational case and adapting an argument used in the parabolic variational inequality setting of \cite{ChipotRodrigues1988} when
 the threshold does not depend on the solution .

 Finally we observe that all the results of this work also apply to linear transport equations
 \begin{equation*}
 \partial_t u+\boldsymbol b\cdot\nabla u+cu=F,
\end{equation*}
under appropriate assumptions on the data $\boldsymbol b$, $c$ and $F$, but since this case presents interesting additional properties
it will be treated in the forthcoming article \cite{RodriguesSantos2012}.

\vspace{5mm}

\section{Assumptions and main results}

Let $\Omega$ be a bounded open subset of $\R^N$ with a Lipschitz boundary $\partial\Omega$ and, for $T>0$,  let $Q=\Omega\times(0,T)$.
We denote by $x=(x_1,\ldots,x_N)$ an arbitrary point of $\Omega$,  $\nabla=(\partial_{ x_1},\ldots,\partial_{ x_N})$
denotes the spatial gradient and $\partial_t$ the partial derivative in order to the variable $t\in\,(0,T)$.
 We represent vectorial functions using a bold symbol, $\boldsymbol{F}=(F_1,\ldots,F_N)$. Its
 divergence is then $\nabla\cdot\boldsymbol{F}=\partial_{x_i}F_i=div\boldsymbol F$.

Assume $\boldsymbol{\Phi}=\boldsymbol{\Phi}(x,t,u):\bar{Q}\times\R\rightarrow\R^N$,
$f=f(x,t,u):\bar{Q}\times\R\rightarrow\R$  are such that, for any $R>0$,
\begin{equation}\label{fif}
\boldsymbol{\Phi}\in\,\boldsymbol W^{2,\infty}\big(Q\times(-R,R)\big),\quad f\in\, W^{1,\infty}\big(Q\times(-R,R)\big),
\end{equation}
and $\nabla\cdot\boldsymbol{\Phi}$ and $f$ are functions with at most linear growth in the variable $u$, uniformly in $(x,t)$, i.e.,
there exist positive constants $c_1$ and $c_2$ such that, for  a.e.   $(x,t)\in Q$ and all $u\in\R$,
\begin{equation}
\label{phif}
\big|\big(\nabla\cdot\boldsymbol{\Phi}\big)(x,t,u)+ f(x,t,u)\big|\le c_1|u|+c_2.
\end{equation}

For the gradient constraint  $G=G(x,u):\Omega\times\R\rightarrow\R$ we assume
\begin{equation}\label{rg}
G\in\, C\big(\R;L^\infty(\Omega)\big),\quad G(x,u)\ge\lambda>0\quad\mbox{ a.e. }x\in\Omega,\ \forall u\in\R,
\end{equation}
in particular $G$ is uniformly continuous in $u$ and bounded in $x$ in each compact of $\R$.

Given a positive function $\varphi\in L^\infty(\Omega)$, we define
\begin{equation*}
\K_{\varphi}=\big\{v\in H^1_0(\Omega):|\nabla v|\le \varphi\ \mbox{ a.e. in }\Omega\big\}.
\end{equation*}

\begin{thm}\label{thm:main}
Under the assumptions (\ref{fif})-(\ref{rg}), for each
\begin{equation}\label{u0}
 u_0\in \K_{G(u_0)}\cap C_0(\bar{\Omega}),
\end{equation}
there exists at least a
 quasivariational solution
 \begin{equation*}
 u\in L^\infty\big(0,T;W^{1,\infty}_0(\Omega)\big)\cap C(\bar{Q}),\qquad  \partial_t u\in L^\infty\big(0,T;M(\Omega)\big),
 \end{equation*}
to the problem
\begin{eqnarray}\label{iqv}
\left\{\begin{array}{l}
u(t)\in\K_{G(u(t))},\  u(0)=u_0,\text{ and }\mbox{ for a.e. }t\in(0,T),\vspace{3mm}\\
\displaystyle{\int_\Omega\partial_tu(t)(v-u(t))+\int_\Omega\boldsymbol{\Phi}(u(t))\cdot\nabla(v-u(t))}\vspace{3mm}\\
\displaystyle{\hspace{6cm}\ge\int_\Omega f(u(t))(v-u(t)),\qquad\forall v\in\K_{G(u(t))}.}
\end{array}
\right.
\end{eqnarray}
\end{thm}

Here $M(\Omega)$ denotes the space of bounded measures in $\Omega$, the first integral in (\ref{iqv}) is interpreted in
the duality  between the spaces $M(\Omega)$ and $C_0(\bar{\Omega})$, which denotes
the space of continuous functions on $\bar\Omega$ which vanish on $
\partial\Omega$. The other integrals are interpreted in the usual Lebesgue sense in $\Omega$.

We shall prove this theorem using the vanishing viscosity approach. In particular, we shall obtain the existence of solution to the following
 parabolic
quasivariational inequality with nonlinear convection and nonlinear reaction terms,  extending a result of \cite{RodriguesSantos2000},
\begin{eqnarray}
\label{iqvdelta}
\left\{\begin{array}{l}
u^\delta(t)\in\K_{G(u^\delta(t))},\ u^\delta(0)=u_0^\delta\mbox{ and for a.e. }t\in(0,T),\vspace{3mm}\\
\displaystyle{\int_\Omega\partial_tu^\delta(t)(v-u^\delta(t))+\delta\int_\Omega
\nabla u^\delta(t)\cdot\nabla(v-u^\delta(t))+\int_\Omega\boldsymbol{\Phi}(u^\delta(t))\cdot\nabla(v-u^\delta(t))}\vspace{3mm}\\
\hspace{5,5cm}\displaystyle{\ge\int_\Omega f(u^\delta(t))(v-u^\delta(t)),\qquad\forall v\in\K_{G(u^\delta(t))}.}
\vspace{1mm}\end{array}
\right.
\end{eqnarray}

\begin{prop} \label{prop:delta}
Under the assumptions (\ref{fif})-(\ref{rg}) and  if $ u_0^\delta\in \K_{G(u_0)}\cap C_0(\bar{\Omega})$ is such that
$\Delta u_0^\delta\in M(\Omega)$, there exists a function
\begin{equation*}u^\delta \in L^\infty\big(0,T;W^{1,\infty}_0(\Omega)\big)\cap C(\bar{Q}),\qquad  \partial_t u^\delta\in L^\infty\big(0,T;M(\Omega)\big),
\end{equation*}
 which
 is solution to the parabolic
quasivariational inequality (\ref{iqvdelta}).
\end{prop}

When the function $G$ does not depend on $u$ the conclusions are stronger and we can obtain the uniqueness of solution in the corresponding problem,
which is then a variational inequality.

\begin{thm}\label{thm:iv}
Under the assumptions (\ref{fif}), (\ref{phif}) and
\begin{equation*}\label{9'}
\tag{8'}
 0<\lambda\le G(x)\le \Lambda\mbox{ for a.e. }x\in\Omega,
\end{equation*}
\begin{equation*}
u_0\in\K_G,
\end{equation*}
there exists a unique function
\begin{equation*}
\vartheta\in L^\infty\big(0,T;W^{1,\infty}_0(\Omega)\big)\cap C(\bar{Q}),\qquad
\partial_t \vartheta\in L^\infty\big(0,T;M(\Omega)\big)\cap L^2(Q),
\end{equation*}
such that $\vartheta$ is the  solution of the first order variational inequality
\begin{eqnarray}
\label{iv}
\left\{\begin{array}{l}
\displaystyle{\vartheta(t)\in \K_G,\ \vartheta(0)=u_0\mbox{ and  for a.e. }t\in(0,T)}, \vspace{3mm}\\
\displaystyle{\int_\Omega\partial_t\vartheta(t)(v-\vartheta(t))+\int_\Omega\boldsymbol{\Phi}(\vartheta(t))\cdot\nabla(v-\vartheta(t))}\vspace{3mm}\\
\hspace{6,5cm}\displaystyle{\ge\int_\Omega f(\vartheta(t))(v-\vartheta(t)),\qquad\forall v\in\K_{G}.}
\vspace{1mm}\end{array}
\right.
\end{eqnarray}
\end{thm}

We prove the stabilization in time (for subsequences $t_n\rightarrow\infty$) to steady-state solutions, both in
the quasivariational and in the variational cases.
Given  functions $f^\infty=f^\infty(x,u):\bar\Omega\times\R\longrightarrow\R$ and $\boldsymbol \Phi=\boldsymbol\Phi(x,u):
\bar\Omega\times\R\longrightarrow\R^N$, we consider the stationary quasivariational inequality: to find
\begin{eqnarray}\label{iqvs}
\left\{\begin{array}{l}
u^\infty\in \K_{G(u^\infty)}\vspace{3mm}\\
\displaystyle{\int_\Omega\boldsymbol{\Phi}(u^\infty)\cdot\nabla(v-u^\infty)\ge\int_\Omega f^\infty(u^\infty)(v-u^\infty),\qquad\forall v\in\K_{G(u^\infty)}.}
\end{array}
\right.
\end{eqnarray}

\begin{thm} \label{asiqv}
 Assume (\ref{fif}), (\ref{phif}), (\ref{u0}), $\boldsymbol\Phi$ is time independent and $f$ is strictly decreasing in $u$, i.e.,
\begin{equation}\label{Phif}
\partial_t\boldsymbol\Phi=0,\qquad \partial_u f\le -\mu<0
\end{equation}
and
\begin{equation}\label{gu}
 0<\lambda\le G(x,u)\le \Lambda\quad\mbox{ for a.e. }x\in\Omega,\quad \mbox{ for all }u\in\R
\end{equation}
or there exists $M>0$ such that, for all $R\ge M$,
\begin{equation}\label{supsubsol}
\big(\nabla\cdot\boldsymbol\Phi\big)(x,R)+f(x,t,R)\le 0,\qquad \big(\nabla\cdot\boldsymbol\Phi\big)(x,-R)+f(x,t,-R)\ge 0.
\end{equation}

Setting $\xi_R(t)=\displaystyle\int_\Omega\sup_{|u|\le R}\big|\partial_tf(x,t,u)\big|\,dx$, we suppose, in addition, that for $R\ge R_0$,
\begin{equation}\label{xiR}
\sup_{0<t<\infty}\int_t^{t+1}\xi_R(\tau)d\tau\le C_R\qquad\text{and}\qquad\int_t^{t+1}\xi_R(\tau)d\tau\underset{t\rightarrow\infty}{\longrightarrow}0,
\end{equation}
being $C_R$ a positive constant and
\begin{equation}\label{fff}
f(x,t,u)\underset{t\rightarrow\infty}{\longrightarrow}f^\infty(x,u),\quad\mbox{ for all }|u|\le R\mbox{ and  a.e. }x\in\Omega.
\end{equation}

Then problem  (\ref{iqvs}) has a solution which is the weak\,-\,$*$ limit in $W^{1,\infty}_0(\Omega)$  and strong limit in
$C^{0,\alpha}(\bar\Omega)$, $0\le\alpha<1$, for $t_n\rightarrow\infty$, of a sequence $\{u(t_n)\}_n$,
being $u$ a global solution of problem  (\ref{iqv}).
\end{thm}

These asymptotic results still hold when $G$ is independent of $u$. However, in this case, when we can guarantee the
 uniqueness of the stationary solution for the variational inequality, we may conclude the convergence of the whole sequence, i.e.,
\begin{equation*}
\vartheta(t)\underset{t\rightarrow\infty}{\lraup}\vartheta^\infty\quad\text{in }W^{1,\infty}(\Omega)\text{ weak\,-\,}*,
\end{equation*}
 and, by compactness, also strongly in spaces of Hölder continuous functions.

\begin{prop}\label{ivs1o}
Suppose that $\boldsymbol\Phi$, $f^\infty$ satisfy  (\ref{Phif}), (\ref{fif}), (\ref{phif}), $G=G(x)$ satisfies (\ref{9'}).

Then the following stationary variational inequality
\begin{eqnarray}\label{ivinfty}
\left\{\begin{array}{l}
\displaystyle{\vartheta^\infty\in \K_G,} \vspace{3mm}\\
\displaystyle{\int_\Omega\boldsymbol{\Phi}(\vartheta^\infty)\cdot\nabla(v-\vartheta^\infty)
\ge\int_\Omega f^\infty(\vartheta^\infty)(v-\vartheta^\infty),\qquad\forall v\in\K_{G},}
\vspace{1mm}\end{array}
\right.
\end{eqnarray}
 has a unique solution $\vartheta^\infty$.
\end{prop}

With additional assumptions we may even estimate the order of convergence of $\vartheta(t)$ to $\vartheta^\infty$.

\begin{thm}\label{asymp}
Under the assumptions of  Proposition \ref{ivs1o}  and
%
\begin{equation*}
\eta_{_M}(t)=\int_t^{t+1}\int_\Omega\sup_{|u|\le M}\big|f(x,\tau,u)-f^\infty(x,u)\big|\,dx\,d\tau\underset{t\rightarrow\infty}{\longrightarrow}0,
\end{equation*}
for  $M\ge M_{\Lambda}$, if $\vartheta$ denotes the solution of the variational inequality (\ref{iv}) and $\vartheta^\infty$ the solution of the variational inequality
 (\ref{iqvs}) then
\begin{equation*}
\vartheta(t)\underset{t\rightarrow\infty}{\longrightarrow}\vartheta^\infty\qquad\text{in }C^{0,\alpha}(\bar\Omega)\ \text{  for }0\le \alpha<1.
\end{equation*}

In addition, if  $\eta_{_M}(t)=O(e^{-\gamma t})$, $\gamma>0$, when $t\rightarrow\infty$,
then, for any $0\le\alpha<1$,
\begin{equation*}
\|\vartheta(t)-\vartheta^\infty\|_{C^{0,\alpha}(\bar\Omega)}=O\big(e^{-\gamma' t}\big),\qquad\text{as }t\rightarrow\infty,
\end{equation*}
where $\gamma'=\big(\frac{1-\alpha}{n+1}\big)\nu$ with $\nu=\min\{\mu,\gamma\}$.
\end{thm}

\vspace{5mm}

\section{The {\em a priori} estimates for the approximating problem}\label{srg}

For $0<\eps<1$, $0<\delta<1$, we consider the approximating quasilinear parabolic problem for $u^{\eps\delta}$ (denoted, for simplicity, also
by $w$),
\begin{eqnarray}
\label{aprox}
\left\{\begin{array}{l}
\partial_tw-\nabla\cdot\left[\delta k_\eps\big(|\nabla w|^2-G_{\eps}^2(w)\big)\nabla w+\boldsymbol{\Phi}(w)\right]=f(w),\qquad\mbox{ in }Q,\vspace{3mm}\\
w(0)=u_0^{\eps\delta}\mbox{ in }\Omega,\qquad w=0\mbox{ on }\partial\Omega\times(0,T),\vspace{1mm}
\end{array}\right.
\end{eqnarray}
where $k_\eps$ is a smooth function such that $k_\eps\big(s)=1$ if $s\le 0$ and
 $k_\eps\big(s)=e^{\frac{s}\eps}$ if $s\ge\eps$, $G_{\eps}=G*\rho_\eps\ge\lambda>0$,
  being $\rho_\eps$ a mollifier in $(x,u)$ and $u_0^{\eps\delta}$ an approximation of $u_0$, belonging to $\mathscr{D}(\Omega)$,
  such that $|\nabla u_0^{\eps\delta}|\le|\nabla u_0|$ and
$\|\Delta u_0^{\eps\delta}\|_{L^1(\Omega)}\le C_0/\delta$, where $C_0$ is a constant independent of $\eps$ and $\delta$.

General parabolic theory for quasilinear non-degenerate equations (\cite{Krylov1987,LSU1968}) yields the existence of a unique solution $
u^{\eps\delta}\in C^2(Q)\cap \,C(\bar{Q})$ to (\ref{aprox}).

We  obtain several {a priori} estimates for the solution $u^{\eps\delta}$ of the approximating problem, independent of $\eps$ and $\delta$.
\begin{prop}\label{cc}
Under the assumptions (\ref{fif})-(\ref{u0}),
the solution $u^{\eps\delta}$
of  problem (\ref{aprox}) is such that,
\begin{enumerate}
\item  for a.e. $(x,t)\in Q$,
\begin{equation}\label{n}
\big|u^{\eps\delta}(x,t)\big|\le M=\inf_{\lambda>b_1}\ e^{\lambda T}\left\{\max_{\bar{\Omega}}|u_0|+1,
\left(\frac{b_2}{\lambda-b_1}\right)^{\frac12}\right\},
\end{equation}
where $b_1=c_1+\frac12$ and $b_2=\frac12\,c_2^2$, $c_1$ and $c_2$ positive constants defined in (\ref{phif}).
In particular, $M$ is independent of $\eps$ and $\delta$;
\item there exist positive constants $C_1, C_2$ and $K$, such that
\begin{equation}
\label{ut}
\big\|\partial_t u^{\eps\delta}\big\|_{L^\infty(0,T;L^1(\Omega))}\le e^{C_2T}\big(C_1+K\|\nabla u^{\eps\delta}\|_{L^1(Q)}\big);
\end{equation}
\item there exists a positive constant $C$  such that
\begin{equation*}
\big\|k_\eps\big(|\nabla u^{\eps\delta}|^2-G_{\eps}^2(u^{\eps\delta})\big)\big\|_{L^\infty(0,T;L^1(\Omega))}\le \frac1{\delta^2}C;
\end{equation*}
\item  for all $2\le q<\infty$, there exists a positive constant $D_q$ such that
\begin{equation}
\label{gd}
\big\|\nabla u^{\eps\delta}\big\|_{L^\infty(0,T;L^q(\Omega))}\le \frac1{\delta^2}D_q,
\end{equation}
where $D_q$ depends on $q$ and is independent of $\eps$ and $\delta$.
\end{enumerate}
\end{prop}

\vspace{3mm}

\noindent{\em Proof of {\em Proposition \ref{cc}}--1.}

\noindent By Theorem 2.9 of \cite{LSU1968}, page 23, if a function $w$ is a classical solution of the equation
\begin{equation}
\label{elsu}
\partial_tw-a_{ij}(x,t,w,\nabla w)\partial_{x_ix_j}w+a(x,t,w,\nabla w)=0,
\end{equation}
with $a_{ij}$ and $a$ bounded and satisfying the following two  conditions,
\begin{equation*}
a_{ij}(x,t,w,0)\xi_i\xi_j\ge 0,\qquad w\,a(x,t,w,0)\ge - b_1w^2-b_2,
\end{equation*}
for all $\xi=(\xi_1,\ldots,\xi_N)\in\R^N$, with $b_1$ and $b_2$ nonnegative constants, then
\begin{equation*}
\max_{Q}\big|w(x,t)\big|\le \inf_{\lambda>b_1}\ e^{\lambda T}
\left\{\max_{\Omega\times\{0\}\cup\,\partial\Omega\times(0,T)}|w|,\left(\frac{b_2}{\lambda-b_1}\right)^{\frac12}\right\}.
\end{equation*}

The equation in (\ref{aprox}) is of the type of (\ref{elsu}), with
\begin{equation*}
a_{ij}(x,w,\nabla w)=2\,\delta k'_\eps\big(|\nabla w|^2-G_{\eps}^2(w)\big)\partial_{x_i}w\,\partial_{x_j}w+
\delta k_\eps\big(|\nabla w|^2-G_{\eps}^2(w)\big)\delta_{ij},
\end{equation*}
\begin{multline*}a(x,t,w,\nabla w)=2\,\delta k'_\eps\big(|\nabla w|^2-G_{\eps}(w)\big)G_{\eps}(w)\big[\partial_uG_{\eps}(w)|\nabla w|^2
+\nabla G_\eps(w)\cdot\nabla w\big]\\
-(\nabla\cdot\boldsymbol\Phi)(x,t,w)-(\partial_u\boldsymbol\Phi)(x,t,w)\cdot\nabla w-f(x,t,w)
\end{multline*}
and
\begin{equation*}
w\,a(x,t,w,0)=-w\,\big((\nabla\cdot\boldsymbol\Phi)(x,t,w)+f(x,t,w)\big).
\end{equation*}

By (\ref{phif}),
$ |(\nabla\cdot\boldsymbol\Phi)(x,t,w)+f(x,t,w)|\le 2\,c_1|w|+2\,c_2$
and so $
wa(x,t,w,0)\ge -b_1w^2-b_2,$
for $b_1=2\,c_1+1$ and $b_2=c_2^2$.
We also have $\displaystyle{|w(x,t)|\le \max_{\bar{\Omega}}|u_0|+1}$ for a.e. $(x,t)\in \Omega\times\{0\}\cup\,\partial\Omega\times(0,T)$
and the conclusion follows. \hfill{$\square$}

\vspace{3mm}

\noindent{\em Proof of {\em Proposition \ref{cc}}--2.}

\noindent Setting $v=\partial_tw$ and differentiating the first equation of (\ref{aprox}) in order to $t$, we get
\begin{multline}
\label{v}
\displaystyle{\partial_tv-\Big[\delta k'_\eps\big(|\nabla w|^2-G^2_{\eps}(w)\big)\big(2\,\partial_{x_i}w\,
\partial_{x_j}w\,\partial_{x_j}v-2\,G_{\eps}(w)
(\partial_uG_{\eps})(w)\,v\,\partial_{x_i}w\big)}\hspace{2cm}\\
\displaystyle{+\delta k_\eps\big(|\nabla w|^2-G^2_{\eps}(w)\big)\partial_{x_i}v+(\partial_t\Phi_i)(w)+(\partial_u\Phi_i)(w)v\Big]_{x_i}=
(\partial_tf)(w)+(\partial_uf)(w)v}.
\end{multline}

We define the signal function by sgn$^0(0)=0,$ sgn$^0(\tau)=1$ if
$\tau>0$ and sgn$^0(\tau)=-1$ if $\tau<0$ and we consider  the approximation by a sequence of $C^1$ increasing functions,
$s_\zeta:\R\rightarrow\R$ such that
\begin{equation}\label{zetaz}
0\le s'_\zeta(\tau)\le\frac{C}\zeta\ \text{ for all }\tau\in\R,\qquad s_{\zeta}(\tau)=
\mbox{sgn}^0(\tau)\ \text{ for all }\tau\in \{0\}\cup\,\R\setminus]-\zeta,\zeta[.
\end{equation}

Setting $\displaystyle{S_\zeta(\tau)
=\int_0^\tau s_\zeta(\sigma)d\sigma}$, we have
\begin{eqnarray}
\label{hdelta}
\lim_{\zeta\rightarrow 0}\tau s'_\zeta(\tau)=0,\quad
\lim_{\zeta\rightarrow 0}S_\zeta(\tau)=|\tau|,\quad \forall \tau\in\R.
\end{eqnarray}

As $v=0$ on $\partial\Omega\times]0,t[$ we also have $s_\zeta(v)=0$ on $\partial\Omega\times]0,t[$. If we multiply
(\ref{v}) by $s_\zeta(v)$ and integrate over $Q_t=\Omega\times(0,t)$, we obtain
\begin{alignat}{1}
\label{zeta}
\nonumber\int_{Q_t}\partial_tv\, s_\zeta(v)&+2\,\delta\int_{Q_t} k'_\eps\big(|\nabla w|^2-G^2_{\eps}(w)\big)
\partial_{x_i}w\,\partial_{x_j}w\,\partial_{x_j}v\,\partial_{x_i}v\,s'_\zeta(v)\\
 &-2\,\delta \int_{Q_t} k'_\eps\big(|\nabla w|^2-G^2_{\eps}(w)\big)\,
G_{\eps}(w)\,(\partial_uG_{\eps})(w)\,
v\,\partial_{x_i}w\,\partial_{x_i}v\,s'_\zeta(v)\\
\nonumber &+\delta \int_{Q_t}k_\eps\big(|\nabla w|^2-G^2_{\eps}(w)\big)\,|\nabla v|^2\,
s'_\zeta(v)
 -\int_{Q_t}\,(\partial_{x_i}\partial_t\Phi_i)(w)\,s_\zeta(v)\\
\nonumber &-\int_{Q_t}(\partial_u\partial_{t}\Phi_i)(w)\, \partial_{x_i}w\,s_\zeta(v)+
\int_{Q_t}\,(\partial_u\boldsymbol{\Phi})(w)\cdot\nabla v\,v\,s'_\zeta(v)\\
\nonumber &=\int_{Q_t}\,(\partial_t f)(w)\, s_\zeta(v)+
\int_{Q_t}\, (\partial_uf)(w)\,v\, s_\zeta(v).
\end{alignat}

We  note that
\begin{equation*}
v(x,0)=\partial_tw(x,0)=\delta\Delta u^{\eps\delta}_0+(\nabla\cdot\boldsymbol{\Phi})(x,0,u^{\eps\delta}_0)+(\partial_u\boldsymbol\Phi)(x,0,u^{\eps\delta}_0)\cdot
\nabla u^{\eps\delta}_0+f(x,0,u^{\eps\delta}_0).
\end{equation*}

Let us analise each  term of the left-hand side of the equality (\ref{zeta}).
\begin{equation*}
\int_{Q_t}\partial_tv\,s_\zeta(v)=\int_{Q_t}\frac{d\ }{dt}S_{\zeta}(v)
=\int_\Omega
\big[S_{\zeta}(v(t))-S_{\zeta}(v(0))\big]=\int_\Omega
S_{\zeta}\big(\partial_tw(t)\big)-\int_\Omega S_{\zeta}\big(\partial_tw(0)\big)
\end{equation*}
and
\begin{multline*}
\lim_{\zeta\rightarrow 0}\int_{Q_t}\partial_tv\,s_\zeta(v)=\int_\Omega\big|\partial_tw(t)\big|\\
-\int_\Omega\big|\delta\Delta u^{\eps\delta}_0+
\left(\nabla\cdot\boldsymbol{\Phi}\right)(x,0,u^{\eps\delta}_0)+(\partial_u\boldsymbol\Phi)(x,0,u^{\eps\delta}_0)\cdot\nabla u^{\eps\delta}_0+f(x,0,u^{\eps\delta}_0)\big|.
\end{multline*}

On the other hand,
\begin{multline*}
\displaystyle{2\,\delta\int_{Q_t}\,k'_\eps\big(|\nabla w|^2-G_{\eps}^2(w)\big)\,
\partial_{x_i}w\,\partial_{x_j}w\,\partial_{x_j}v\,\partial_{x_i}v\,s'_\zeta(v)=}\\
\displaystyle{2\,\delta\,\int_{Q_t}\,k'_\eps\big(|\nabla w|^2-G_{\eps}^2(w)\big)\,
(\nabla w\cdot\nabla v)^2\,s'_\zeta(v)\,\geq 0,}
\end{multline*}
and
\begin{equation*}
\delta\,\int_{Q_t}\,k_\eps\big(|\nabla w|^2-G_{\eps}^2(w)\big)\,|\nabla v|^2\,s'_\zeta(v)\,\geq 0,
\end{equation*}
\begin{equation*}
\Big|\int_{Q_t}\,(\partial_t\partial_{x_i}\Phi_i)(w)\,s_\zeta(v)\Big|\,\le K,
\qquad
\Big|\int_{Q_t}\,(\partial_t\partial_u\Phi_i)(w)\,\partial_{x_i}w\,s_\zeta(v)\Big|\,\le\, K\,\|\nabla w\|_{L^1(Q)},
\end{equation*}
where, for $M$ defined in (\ref{n}), $ K=\|\boldsymbol{\Phi}\|_{W^{2,\infty}(Q\times(-M,M))}$
is independent of $\eps$ and $\delta$.
By (\ref{hdelta}) and Lebesgue Theorem
\begin{equation*}
-2\,\delta\,\int_{Q_t}\,k'_\eps\big(|\nabla w|^2-G_{\eps}^2(w)\big)\,G_{\eps}(w)\,
(\partial_uG_{\eps})(w)\,\partial_{x_i}w\,\partial_{x_i}v\,v\,s'_\zeta(v)\,\longrightarrow 0,\mbox{ when
}\zeta\rightarrow 0
\end{equation*}
and
\begin{equation*}
\int_{Q_t}\,(\partial_u\boldsymbol{\Phi})(w)\cdot\nabla v\,v\,s'_\zeta(v)\,\longrightarrow 0,\mbox{ when
}\zeta\rightarrow 0.
\end{equation*}

Gathering all the information above,
\begin{multline}
\label{vvv}
\int_\Omega\,\big|v(t)\big|\le\int_\Omega\,\big|\delta\,\Delta u^{\eps\delta}_0+(\nabla\cdot\boldsymbol{\Phi})(x,0,u^{\eps\delta}_0)
+(\partial_u\boldsymbol\Phi)(x,0,u^{\eps\delta}_0)\cdot\nabla u^{\eps\delta}_0+f(x,0,u^{\eps\delta}_0)\big|\\
+K\big\|\nabla w\big\|_{L^1(Q)}+
\int_{Q_t}\,\big| \partial_tf(w)\big|+\int_{Q_t}\,\big| \partial_uf(w)v\big|,
\end{multline}
and we obtain
\begin{equation*}
\int_\Omega\big|v(t)\big|\le
C_1+K\big\|\nabla w\big\|_{L^1(Q)}+C_2\int_0^t\int_\Omega\big| v\big|.
\end{equation*}

Here, the constant $C_1$ is dependent on $\|u_0\|_{H^1_0(\Omega)}$, $\delta\|\Delta u_0^{\eps\delta}\|_{L^1(\Omega)}$ (bounded by $C_0$),
$\|\boldsymbol{\Phi}\|_{\boldsymbol W^{1,\infty}(Q\times(-M,M))}$,
$\|f\|_{W^{1,\infty}(Q\times(-M,M))}$, the constant $C_2$ depends on $\|f\|_{W^{1,\infty}(Q\times(-M,M)}$
 and $ K$ depends only on $\|\boldsymbol{\Phi}\|_{W^{2,\infty}(Q\times(-M,M))}$.
The conclusion follows immediately by applying the Gronwall inequality.
\hfill{$\square$}

\vspace{3mm}

\begin{rem}
When $\boldsymbol{\Phi}$ is independent of $t$, we may take  $K=0$ in the last theorem.
\end{rem}

\vspace{3mm}

\noindent{\em Proof of {\em Proposition \ref{cc}}--3.}

\noindent We multiply the equation (\ref{aprox}) by $w$ and integrate over  $\Omega$, obtaining, for a.e. $t\in (0,T)$,
\begin{equation*}
\delta\int_\Omega k_\eps(|\nabla w(t)|^2-G_{\eps}^2(w(t)))\,\big|\nabla w(t)\big|^2=\int_\Omega\big(f(w(t))-
\partial_t w(t)\big)\,w(t)-\int_\Omega\boldsymbol\Phi(w(t))\cdot\nabla w(t).
\end{equation*}
So we get for $w=w(t)$,
\begin{multline*}
\delta\int_\Omega\big|\nabla w\big|^2\le\delta\int_\Omega k_\eps(|\nabla w |^2-G_{\eps}^2(w ))\,\big|\nabla w \big|^2\\
\le M\big\|f-\partial_t w\big\|_{L^\infty(0,T;L^1(\Omega))}+\frac1{2\delta}\big\|\boldsymbol\Phi(w)\big\|^2_{L^\infty(0,T;L^2(\Omega))}+
\frac\delta{2}\int_\Omega\big|\nabla w\big|^2,
\end{multline*}
from which we conclude
\begin{equation}\label{2728}
\int_\Omega\big|\nabla w(t)\big|^2\le \int_\Omega k_\eps(|\nabla w(t)|^2-G_{\eps}^2(w(t)))\,\big|\nabla w(t)\big|^2\le \frac{C}{\delta^2},\qquad\text{for a.e. }t\in(0,T).
\end{equation}

As $k_\eps(s)s\ge 0$ for $s\ge 0$, $k_\eps(s)=1$ if $s\le 0$ and $G_\eps\ge\lambda$,
 we obtain the following estimate, for a.e. $t\in(0,T)$,
 \begin{align*}
\lambda^2\int_\Omega k_\eps(|\nabla w(t)|^2&-G_{\eps}^2(w(t)))\le
\int_\Omega k_\eps(|\nabla w(t)|^2-G_{\eps}^2(w(t)))\,G_{\eps}^2(w(t))\\
&=\int_\Omega k_\eps(|\nabla w(t)|^2-G_{\eps}^2(w(t)))\,\big(G_{\eps}^2(w(t))-|\nabla w(t)|^2\big)\\
&\ \ \ \ +
\int_\Omega k_\eps(|\nabla w(t)|^2-G_{\eps}^2(w(t)))\,\big|\nabla w(t)\big|^2\\
&\le \int_{\{|\nabla w(t)|<G_\eps^2(w(t))\}} G_{\eps}^2(w(t))+\frac1{\delta^2}C\le \big(\big\|G\big\|^2_{L^\infty}+1\big)\,\big|\Omega\big|
+\frac1{\delta^2}C.
 \end{align*}

\hfill{$\square$}

\vspace{3mm}

\noindent{\em Proof of {\em Proposition \ref{cc}}--4.}

\noindent Calling
 \begin{equation*}
 A_\eps(t)=\{x\in\Omega:|\nabla w(t)|^2>G_{\eps}^2(w(t))+\eps\}.
 \end{equation*}
we decompose $\Omega=A_\eps(t)\cup\big(\Omega\setminus A_\eps(t)\big)$.
Then
\begin{equation}
\label{leq}
\int_{\Omega\setminus A_\eps(t)}\,\big|\nabla w(t)\,\big|^q
\le\int_\Omega
\big(G_{\eps}^2(w(t)))+\eps\big)^{\frac{q}2}\le \big(\|G\big\|^2_{L^\infty}+1\big)^{\frac{q}2}\,\big|\Omega\big|.
\end{equation}

Recalling (\ref{2728}) and, as  $k_\eps(s)=e^{\frac{s}\eps}\ge \frac{s^q}{q!\,\eps^q}$, for $s>\eps$, we have, for $q$ even integer greater or equal to
 $4$,
\begin{multline}
\label{ge}
\int_{A_\eps(t)}
\frac{\left(|\nabla w(t)|^2 -G_{\eps}^2(w(t)\right))^{\frac{q-2}2}}{\eps^{\frac{q-2}2} \big(\frac{q-2}2\big)!}\,\big|\nabla w(t)\big|^2\\
\le\int_{A_\eps(t)}k_\eps(|\nabla w(t)|^2-G_{\eps}^2(w(t)))\,\big|\nabla w(t)\big|^2
\le \frac{C}{\delta^2}.
\end{multline}

Then, using (\ref{ge})  we find
\begin{align}
\label{wwww}
\nonumber\int_{A_\eps(t)}\big|\nabla &w(t)\big|^q=\int_{A_\eps(t)}\Big[\Big(\big|\nabla w(t)\big|^2-G^2_{\eps}(w(t))\Big)+
G^2_{\eps}(w(t))\Big]^{\frac{q-2}2}\big|\nabla w(t)\big|^2
\\
\nonumber&\le\int_{A_\eps(t)} 2^{\frac{q-2}2-1} \left(\Big(\big|\nabla w(t)\big|^2-G^2_{\eps}(w(t))
\Big)^{\frac{q-2}{2}}+G^{q-2}_{\eps}(w(t))\right)
\big|\nabla w(t)\big|^2\\
&\le 2^{\frac{q-4}2} \eps^{\frac{q-2}2}\,\big(\tfrac{q-2}2\big)!\, \int_{A_\eps(t)}k_\eps\big(|\nabla w(t)|^2-G^2_{\eps}(w(t))\big)\,\big|\nabla w(t)\big|^2\\
\nonumber&\ \ \ \ +
2^{\frac{q-4}2} \int_{A_\eps(t)} G^{q-2}_{\eps}(w(t))\,\big|\nabla w(t)\big|^2\\
\nonumber&\le \frac1{\delta^2}\,C_q+
2^{\frac{q-4}2}\big\|G_{\eps}\big\|_{L^\infty}^{q-2}\big\|\nabla w\big\|^2_{L^\infty(0,T;L^2(\Omega))}.
\end{align}

Using again (\ref{2728}), (\ref{leq}) and (\ref{wwww}), we obtain that
\begin{equation*}\big\|\nabla w\big\|^q_{L^\infty(0,T;L^q(\Omega))}\le\frac{D_q}{\delta^2}.
\end{equation*}
\hfill{$\square$}

\vspace{5mm}

\section{Existence of solutions by letting $\eps\rightarrow 0$ and $\delta\rightarrow 0$}

In this section we prove Proposition \ref{prop:delta} and Theorem \ref{thm:main}. We start with a simple lemma:
let a sequence $\{h^\eps\}$ of functions
belong to $L^\infty(\Omega)$ and converge to $h$, where $h$ and $h^\eps$ are greater or equal to $\lambda>0$ and  define
$\alpha^\eps=\|h-h^\eps\|_{L^\infty(\Omega)}$. Then, for any $v\in\K_h$,
$v^\eps=\frac\lambda{\lambda+\alpha^\eps}v\in\K_{h^\eps}$ and $v^\eps\underset{\eps\rightarrow 0}{\longrightarrow}v$ in $W^{1,\infty}_0(\Omega)$. This strong
approximation of
a given function belonging to $\K_h$ by a sequence of functions belonging to $\K_{h^\eps}$ will be a key tool in the proof of existence of solution to the
 problems (\ref{iqvdelta}) and (\ref{iqv}).

 \vspace{3mm}

\noindent{\em Proof of {\em Proposition \ref{prop:delta}}.}

\noindent By  the estimates (\ref{ut}) and (\ref{gd}), we have
the
uniform (independent of $\eps$) boundedness of $\{\partial_t u^{\eps\delta}\}_\eps$ in $L^\infty(0,T;L^1(\Omega))$ and of
$\{ \nabla u^{\eps\delta} \}_\eps$ in $L^\infty(0,T;L^q(\Omega))$. This implies the convergence, at least for a subsequence,
\begin{equation*}
\partial_tu^{\eps\delta}\lraup\partial_tu^\delta\qquad\mbox{ in }L^\infty(0,T;M(\Omega))\mbox{ weak\,-\,}*\mbox{ when }\eps\rightarrow 0
\end{equation*}
\begin{equation*}
\nabla u^{\eps\delta}\lraup\nabla u^\delta\qquad\mbox{ in }L^\infty(0,T;L^q(\Omega))\mbox{ weak\,-\,}*\mbox{  when }\eps\rightarrow 0.
\end{equation*}

By the compactness of the imbedding $W^{1,q}_0(\Omega)\subset C(\bar{\Omega})$, for $q>N$, and a theorem on compactness of functions
with values in an intermediate  Banach space (see \cite{Simon1987}, page 84)  $\{u^{\eps\delta}\}$ is
relatively compact in $C([0,T];C(\bar{\Omega}))$. So, there exists a subsequence  such that
\begin{equation*}
%
u^{\eps\delta}\longrightarrow u^\delta\qquad \mbox{ in }\qquad C(\bar{Q}),\mbox{ when }\eps\rightarrow 0
\end{equation*}
and, consequently, $G_{\eps}(u^{\eps\delta})\underset{\eps\rightarrow 0}{\longrightarrow}G(u^\delta)$ in $L^\infty(Q)$ and $u^\delta(0)=u_0^\delta$
in $\Omega$.

We are going to prove that $u^\delta(t)\in\K_{G(u^\delta(t))}$, for a.e. $t\in(0,T)$. Set
$A_\eps=\{(x,t)\in Q: |\nabla u^{\eps\delta}(x,t)|\ge G_{\eps}(u^{\eps\delta})+\sqrt\eps\}$.
 Then, by the definitions of $k_\eps$ and  the uniform boundedness (in $\eps$) of
$\{k_\eps\big(|\nabla u^{\eps\delta}|^2-G_{\eps}^2(u^{\eps\delta})\big)\}_\eps$ in $L^1(Q)$,
\begin{equation*}
\big|A_\eps\big|=\int_{A_\eps} 1\le\int_{A_\eps}\frac{k_\eps\big(|\nabla u^{\eps\delta}|^2-G_\eps^2(u^{\eps\delta})\big)}{e^{\frac1{\sqrt{\eps}}}}\le Ce^{-\frac1{\sqrt{\eps}}}
\end{equation*}
and this implies that $\big|A_\eps\big|\underset{\eps\rightarrow 0}{\longrightarrow}0$. So
\begin{align*}
\int_Q\big(|\nabla u^\delta|-G(u^\delta)\big)^+&=\int_Q\liminf_\eps\big(|\nabla u^{\eps\delta}|-G(u^{\eps\delta})-\sqrt\eps\big)^+\\
&\le \liminf_{\eps\rightarrow 0}\int_Q\big(|\nabla u^{\eps\delta}|-G(u^{\eps\delta})-\sqrt\eps\big)\,\cchi_{A_\eps}\\
&\le\lim_{\eps\rightarrow 0}\big\||\nabla u^{\eps\delta}|-G(u^{\eps\delta})\big\|_{L^2(Q)}\,\big|A_\eps\big|^{\frac12}=0
\end{align*}
and  $\big|\nabla u^\delta\big|\le G(u^\delta)$ for a.e. in $Q$, i.e., $u^\delta(t)\in\K_{G(u^\delta(t)}$, for a.e. $t\in(0,T)$.

As we stated at the beginning of this section, given $v^\delta\in  L^\infty(0,T;C(\bar\Omega))$ such
that $v^\delta(t)\in\,\K_{G(u^\delta(t))}$
for a.e. $t\in(0,T)$, we may find
$v^{\eps\delta}\in  L^\infty(0,T;C(\bar\Omega))$ such that $v^{\eps\delta}(t)\in\,\K_{G_\eps(u^{\eps\delta}(t))}$
 and $v^{\eps\delta}\underset{\eps\rightarrow 0}{\longrightarrow} v^\delta$ in $L^\infty(0,T;W^{1,\infty}_0(\Omega)),$
 since $G_\eps(u^{\eps\delta}(t))\underset{\eps\rightarrow0}{\longrightarrow}G(u^\delta(t))$ in $L^\infty(\Omega)$.

Multiplying the first equation of (\ref{aprox}) by $v^{\eps\delta}-u^\delta$ and integrating over $\Omega\times(s,t)$, with $0<s<t<T$, we obtain
\begin{multline*}
\int_s^t\int_\Omega\partial_t u^{\eps\delta} \,\big(v^{\eps\delta} -u^\delta \big)
+\delta\int_s^t\int_\Omega k_\eps\big(|\nabla u^{\eps\delta} |^2-G_{\eps}^2(u^{\eps\delta })\big)
\,\nabla u^{\eps\delta} \cdot\nabla\big(v^{\eps\delta} -u^\delta \big)\\
+\int_s^t\int_\Omega\boldsymbol{\Phi}\big(u^{\eps\delta} \big)\cdot\nabla\big(v^{\eps\delta} -u^\delta \big)
=\int_s^t\int_\Omega f\big(u^{\eps\delta} \big)\,\big(v^{\eps\delta} -u^\delta \big).
\end{multline*}

By the monotonicity of $k_\eps$ and since $v^{\eps\delta}(\tau)\in\K_{G_\eps(u^{\eps\delta}(\tau))}$ for a.e. $\tau\in(0,T)$, we have
\begin{equation*}
\int_s^t\int_\Omega\nabla v^{\eps\delta} \cdot\nabla\big(v^{\eps\delta} -u^\delta \big)
\ge \int_s^t\int_\Omega k_\eps\big(|\nabla u^{\eps\delta} |^2-G_{\eps}^2(u^{\eps\delta} )\big)
\nabla u^{\eps\delta} \cdot\nabla\big(v^{\eps\delta} -u^\delta \big),
\end{equation*}
and  the limit when $\eps\rightarrow 0$
yields
\begin{multline*}
\int_s^t\int_\Omega \partial_t u^\delta \,\big(v^\delta -u^\delta \big)+
\delta\int_s^t\int_\Omega\nabla v^{\delta} \cdot\nabla\big(v^\delta -u^\delta \big)\\
+\int_s^t\int_\Omega \boldsymbol{\Phi}\big(u^\delta \big)\cdot\big(v^\delta -u^\delta \big)
\ge \int_s^t\int_\Omega f\big(u^\delta \big)\,\big(v^{\eps\delta} -u^\delta \big).
\end{multline*}

Using Minty's Lemma and since $s$ and $t$ are arbitrary, we obtain as in \cite{RodriguesSantos2000},
\begin{multline*}
\int_\Omega \partial_t u^\delta(t) \,\big(v^\delta(t) -u^\delta(t) \big)+
\delta\int_\Omega\nabla u^{\delta}(t) \cdot\nabla\big(v^\delta(t) -u^\delta(t) \big)\\
+\int_\Omega \boldsymbol{\Phi}\big(u^\delta(t) \big)\cdot\big(v^\delta(t) -u^\delta(t) \big)
\ge \int_\Omega f\big(u^\delta(t) \big)\,\big(v^{\delta}(t) -u^\delta(t) \big),
\end{multline*}
 for a.e.
 $t\in(0,T)$, which concludes the proof.
\hfill{$\square$}

\vspace{3mm}

\noindent{\em Proof of {\em Theorem \ref{thm:main}.}}

\noindent  The solution $u^\delta$ of problem (\ref{iqvdelta}) satisfies $|\nabla u^\delta|\le G(u^\delta)\le G^*$  a.e. in $Q$ and therefore, by assumption
(\ref{rg}) and estimate (\ref{n}), has a gradient uniformly bounded in $\delta$, i.e.,
\begin{equation}
\label{24}
\big\|u^\delta\big\|_{L^\infty(0,T;W^{1,\infty}_0(\Omega))}\le C
\end{equation}
where $C$ is a constant independent of $\delta$.

Letting $\eps\rightarrow 0$ in the estimate (\ref{ut}), and using (\ref{24}), we obtain also
\begin{equation*}
\big\|\partial_t u^{\delta}\big\|_{L^\infty(0,T;M(\Omega))}\le e^{C_2T}\Big(C_1+K\big\|\nabla u^{\delta}\big\|_{L^1(Q)}\Big)
\le C_3.
\end{equation*}

So, there exists a function $u$ such that, when $\delta\rightarrow 0$,
\begin{equation*}
\partial_tu^\delta\lraup \partial_t u \ \mbox{ in }L^\infty(0,T;M(\Omega))
\mbox{ weak\,-\,}*,\qquad\nabla u^\delta\lraup\nabla u \ \mbox{ in }L^\infty(Q)\mbox{ weak\,-\,}*,
\end{equation*}
and, again by the same compactness argument used in the previous proof, we also have the strong convergence
\begin{equation*}
u^\delta\underset{\delta\rightarrow 0}\longrightarrow u \qquad\mbox{ in }C(\bar{Q}),
\end{equation*}
which, in particular, yields $u(0)=u_0$ in $\Omega$.

Then, as before, for a fixed $t\in(0,T)$, given $v\in L^\infty(0,T;C(\bar\Omega))$ such that $v(t)\in \K_{G(u(t))}$,
 we can find $v^\delta(t)\in L^\infty(0,T;C(\bar\Omega))$ such that $v^\delta(t)\in\K_{G(u^\delta(t))}$
 and $v^\delta\underset{\delta\rightarrow 0}{\longrightarrow}v$
 in $L^\infty(0,T;W^{1,\infty}_0(\Omega))$ and we can take the
 limit  $\delta\rightarrow 0$ in
\begin{multline*}
\displaystyle{\int_s^t\int_\Omega\partial_tu^\delta \,\big(v^\delta-u^\delta \big)+\delta\int_s^t\int_\Omega
\nabla u^\delta \cdot\nabla\big(v^\delta-u^\delta \big)+\int_s^t\int_\Omega\boldsymbol{\Phi}\big(u^\delta \big)\cdot\nabla\big(v^\delta-u^\delta \big)}\\
\ge \int_s^t\int_\Omega f\big(u^\delta \big)\,\big(v^\delta -u^\delta \big).
\end{multline*}

As in previous step, by the arbitrariness of $0<s<t<T$, we conclude, for a.e. $t\in(0,T)$,
\begin{equation*}
\int_\Omega \partial_t u(t) \,\big(v(t) -u(t) \big)+
\int_\Omega \boldsymbol{\Phi}\big(u(t) \big)\cdot\big(v(t) -u(t) \big)
\ge \int_\Omega f\big(u(t) \big)\,\big(v(t) -u(t) \big),\quad\text{a.e. }t\in(0,T).
\end{equation*}

Finally, given any measurable set $\omega\subset Q$, we have
\begin{equation*}
\int_\omega\big|\nabla u\big|\le\liminf_{\delta\rightarrow 0}\int_\omega\big|\nabla u^\delta\big|\le\lim_{\delta\rightarrow 0}
\int_\omega G\big(u^\delta\big)=\int_\omega G\big(u\big)
\end{equation*}
and consequently, for a.e. $t\in(0,T)$ we have
\begin{equation*}
\big|\nabla u(t)\big|\le G\big(u(t)\big)\qquad\mbox{ a.e. in }\Omega,
\end{equation*}
which means that $u(t)\in\K_{G(u(t))}$ for  $t\in(0,T)$, concluding the proof that $u$ solves (\ref{iqv}).
\hfill{$\square$}

\vspace{5mm}

\section{The variational inequality}

When the gradient constraint $G$ is a function independent of $u$, the first order quasivariational inequality becomes a variational inequality. As a consequence,
the solution is unique and more regular.

\vspace{3mm}

\noindent {\em Proof of {Theorem \ref{thm:iv}.}}

\noindent We only need to prove that $\partial_t\vartheta\in L^2(Q)$ and the uniqueness of solution.
Multiply the first equation of problem (\ref{aprox})
by $\partial_tw$ and integrate in $Q$. Then
\begin{equation}\label{dt}
\int_Q\big|\partial_tw\big|^2+\delta\int_Qk_\eps\big(|\nabla w|^2-G_\eps^2\big)\nabla w\cdot\nabla\partial_t w
+\int_Q\boldsymbol\Phi(w)\cdot\nabla \partial_tw=\int_Qf(w)\,\partial_tw.
\end{equation}

We remark that
\begin{equation*}
\int_Q\boldsymbol\Phi(w)\cdot\nabla \partial_tw=-\int_Q\big(\nabla\cdot\boldsymbol\Phi\big)(w)\,\partial_tw-
\int_Q\partial_u\boldsymbol\Phi(w)\cdot\nabla w\,\partial_tw
\end{equation*}
and
we obtain
\begin{equation*}
\Big|\int_Q \Big(f(w)+\big(\nabla\cdot\boldsymbol\Phi\big)(w)+\partial_u\boldsymbol\Phi(w)\cdot\nabla w\Big)\partial_tw\Big|\le\frac12\int_Q\big|
\partial_tw\big|^2+C'\|\nabla w\|_{L^2(Q)},
\end{equation*}
where $C'$ is a constant depending only on $M$ and on the assumptions on $\boldsymbol\Phi$ and $f$.

Let $K_\eps(s)=\displaystyle\int_0^sk_\eps(\tau)d\tau$. We have
$K_\eps\big(\big|\nabla w(0)\big|^2-G_\eps^2\big)=\big|\nabla u^{\delta}_0\big|^2-G_\eps^2$ and also

\begin{equation*}
-K_\eps\big(\big|\nabla w(t)\big|^2-G_\eps^2\big)\le G_\eps^2\qquad\text{ a.e. in }Q,
\end{equation*}
since $K_\eps\big(\big|\nabla w(t)\big|^2-G_\eps^2\big)=\big|\nabla w(t)\big|^2-G_\eps^2\ge -G_\eps^2$ if $\big|\nabla w(x,t)\big|^2\le G_\eps^2(x)$
and otherwise $K_\eps\big(\big|\nabla w(t)\big|^2-G_\eps^2\big)\ge 0$.

>From (\ref{dt}) we then obtain,
\begin{equation}\label{acima}
\int_Q\big|\partial_t u^{\eps\delta}\big|^2\le2C'\|\nabla u^{\eps\delta}\|_{L^2(Q)}+\delta C''.
\end{equation}

Letting first $\eps\rightarrow 0$ in (\ref{acima}), we have
\begin{equation*}
\int_Q\big|\partial_t u^{\delta}\big|^2\le 2C'\|\nabla u^{\delta}\|_{L^2(Q)}+\delta C''\le 2C'\|G\|_{L^2(Q)}+\delta C''.
\end{equation*}

Finally, letting $\delta$ tend to zero,
we conclude  that $\partial_t\vartheta\in L^2(Q)$.

To prove the uniqueness of solution we consider  a sequence $s_\zeta:\R\rightarrow\R$
 of $C^1$ increasing odd functions, approximating pointwise the function sgn$^0$, defined as in (\ref{zetaz}). If $\vartheta_1$ and $\vartheta_2$ are two solutions of the variational inequality
 (\ref{iv}) we define $v_1=\vartheta_1+\zeta^2 s_\zeta(\vartheta_2-\vartheta_1)$, for $\zeta$ small. As
\begin{align*}
\big|\nabla v_1\big|&=\big|\nabla \vartheta_1+\zeta^2 s'_\zeta\big(\vartheta_2-\vartheta_1\big)\nabla\big(\vartheta_2-\vartheta_1\big)\big|\\
&=\big|\big(1-\zeta^2s'_\zeta(\vartheta_2-\vartheta_1)\big)\nabla \vartheta_1+\zeta^2s'_\zeta(\vartheta_2-\vartheta_1)\nabla \vartheta_2\big|\\
&\le \big(1-\zeta^2s'_\zeta(\vartheta_2-\vartheta_1)\big)\,G+\zeta^2s'_\zeta(\vartheta_2-\vartheta_1)\,G=G,
\end{align*}
since $1-\zeta^2s'_\zeta(\vartheta_2-\vartheta_1)>0$, we have $v_1\in\K_G$.

Substituting $v_1$ in the variational inequality (\ref{iv}) satisfied by $\vartheta_1$, we get
\begin{multline*}
\int_\Omega\partial_t\vartheta_1(t)\,s_\zeta\big(\vartheta_2(t)-\vartheta_1(t)\big)+
\int_\Omega s'_\zeta\big(\vartheta_2(t)-\vartheta_1(t)\big)\,\boldsymbol\Phi\big(\vartheta_1(t)\big)\cdot\nabla\big(\vartheta_2(t)-\vartheta_1(t)\big)
\\
\ge \int_\Omega f\big(\vartheta_1(t)\big)\,s_\zeta\big(\vartheta_2(t)-\vartheta_1(t)\big).
\end{multline*}

Using now $v_2=\vartheta_2+\zeta^2 s_\zeta(\vartheta_1-\vartheta_2)$ as test function in the problem for $\vartheta_2$, we obtain the following
inequality for
$\bar u=\vartheta_1-\vartheta_2$,
\begin{multline}\label{44}
\int_\Omega \partial_t\bar u(t)\,s_\delta(u(t))+\int_\Omega s'_\zeta\big(\bar u(t)\big)\, (\boldsymbol\Phi\big(\vartheta_1(t)\big)
-\boldsymbol\Phi\big(\vartheta_2(t))\big)\cdot\,\nabla\bar u(t)\\
\le\int_\Omega\big(f(\vartheta_1(t))-f(\vartheta_2(t))\big)\,s_\zeta\big(\bar u(t)\big).
\end{multline}

By assumption (\ref{fif}) there exists $L_{\boldsymbol\Phi}>0$ such that
$|\boldsymbol\Phi(x,u)-\boldsymbol\Phi(x,v)\le L_{\boldsymbol\Phi}|u-v|$. So,
 by (\ref{hdelta}),
\begin{equation*}
\Big|\int_\Omega s'_\zeta(\bar u(t)) \big(\boldsymbol\Phi(\vartheta_1(t))-\boldsymbol\Phi(\vartheta_2(t))\big)\cdot\nabla\bar u(t)\Big|
\le L_{\boldsymbol\Phi}
\int_\Omega\big|\nabla u(t)\big|\,s'_\zeta(\bar u(t))
\,\big|\bar u(t)\big|\underset{\zeta\rightarrow 0}{\longrightarrow}0,
\end{equation*}
applying the dominated convergence theorem and  $\displaystyle\lim_{\tau\rightarrow 0}\tau s'_\zeta(\tau)=0$.
For some positive constant $L_f>0$, we also have
$|f(v_1)-f(v_2)|\le L_f|v_1-v_2|$,  and so
\begin{equation*}
\Big|\int_\Omega\big(f(\vartheta_1(t))-f(\vartheta_2(t))\big)\,s_\zeta(\bar u(t))\Big|\le L_f\int_\Omega\big|\bar u(t)\big|.
\end{equation*}

Then,  integrating (\ref{44}) between $0$ and $t$ and letting $\zeta\rightarrow 0$, we obtain, recalling (\ref{hdelta}) and $\bar u(0)=0$
\begin{equation*}
\lim_{\zeta\rightarrow 0}\int_\Omega S_\zeta(\bar u(t))=\int_\Omega\big|\bar u(t)\big|\le L_f\int_0^t\int_\Omega\big|\bar u\big|.
\end{equation*}

The uniqueness follows by the Gronwall inequality.
\hfill{$\square$}

\vspace{3mm}

\noindent {\em Proof of {\em Proposition \ref{ivs1o}.}}

\noindent The existence of solution to (\ref{ivinfty}) is a consequence of Theorem \ref{asiqv}, since the assumption (\ref{gu}) reduces to (\ref{9'})
when $G$ does not depend on $u$.

The uniqueness can be shown as in the evolutive case: suppose that $\vartheta_1$ and $\vartheta_2$ are two solutions of (\ref{ivinfty}); remarking that
if $s_\zeta$ denotes a $C^1$ approximation of the sgn$^0$ function as in (\ref{zetaz}), we have, in $L^p(\Omega)$, $1\le p<\infty$,
\begin{equation*}
s_\zeta(\vartheta_1-\vartheta_2)\underset{\zeta\rightarrow 0}{\longrightarrow}\text{sgn}^0(\vartheta_1-\vartheta_2)=\left\{\begin{array}{ll}
1&\ \text{ in }\{\vartheta_1>\vartheta_2\},\vspace{3mm}\\
0&\ \text{ in }\{\vartheta_1=\vartheta_2\},\vspace{3mm}\\
-1&\ \text{ in }\{\vartheta_1<\vartheta_2\}.
\end{array}
\right.
\end{equation*}

Using $v=\vartheta_1-\zeta^2s_\zeta(\vartheta_1-\vartheta_2)$ in (\ref{ivinfty}) for $\vartheta^\infty=\vartheta_1$ and $v=\vartheta_2+\zeta^2
s_\zeta(\vartheta_1-\vartheta_2)$ in (\ref{ivinfty}) for $\vartheta^\infty=\theta_2$, we find, recalling that $\boldsymbol\Phi$ is Lipschitz in $u$,
\begin{equation*}
\int_\Omega\Big(f^\infty(\vartheta_2)-f^\infty(\vartheta_1)\Big)s_\zeta(\vartheta_1-\vartheta_2)\le
2L_{\boldsymbol\Phi}\int_\Omega s'_\zeta(\vartheta_1-\vartheta_2)\,\big|\vartheta_1-\vartheta_2\big|\, G.
\end{equation*}

Taking the limit $\zeta\rightarrow 0$, we obtain
\begin{equation*}
\int_{\{\vartheta_1>\vartheta_2)\}}\Big(f^\infty(\vartheta_2)-f^\infty(\vartheta_1)\Big)-
\int_{\{\vartheta_1<\vartheta_2)\}}\Big(f^\infty(\vartheta_2)-f^\infty(\vartheta_1)\Big)\le 0
\end{equation*}
and, since $f^\infty$ is strictly decreasing,  we conclude that
\begin{equation*}
\big|\{\vartheta_1>\vartheta_2)\}\big|=
\big|\{\vartheta_1<\vartheta_2)\}\big|=0
\end{equation*}
 and so $\vartheta_1=\vartheta_2$ a.e. in $\Omega$.\hfill{$\square$}.

\vspace{5mm}

\section{The asymptotic behavior in time}

{\em Proof of {Theorem \ref{asiqv}.}}

\noindent We consider again the approximating problem (\ref{aprox}) defined in $\Omega\times(0,\infty)$ in order to show that
$\partial_t u^{\eps\delta}(t)$ vanishes uniformly in $\eps$ and $\delta$ in $M(\Omega)$ as $t\rightarrow\infty$.

We multiply (\ref{v}) by $e^{\mu t}s_\zeta(v(t))$, with $s_\zeta$ defined in
(\ref{zetaz}) and satisfying (\ref{hdelta}) and $v=\partial_t u^{\eps\delta}(t)$,  using integration by parts in $\Omega\times(s,t)$, $0\le s<t$:
\begin{equation*}
\int_s^t\int_\Omega e^{\mu \tau}\partial_tv(\tau)s_\zeta(v(\tau))\,d\tau=e^{\mu t}\int_\Omega S_\zeta(v(t))- e^{\mu s}\int_{\Omega}S_\zeta(v(s))
-\int_s^t\int_\Omega\mu e^{\mu\tau}S_\zeta(v(\tau))d\tau.
\end{equation*}

Similarly to (\ref{vvv}), where now $K=0$ since $\boldsymbol\Phi$ is independent of $t$ by (\ref{Phif}), after taking $\zeta\rightarrow 0$ we obtain
\begin{multline}\label{37}
e^{\mu t}\int_\Omega \big|\partial_t u^{\eps\delta}(t)\big|\le
e^{\mu s}\int_\Omega \big|\partial_t u^{\eps\delta}(s)\big|+\int_s^t\int_\Omega\mu e^{\mu\tau}\big|\partial_t u^{\eps\delta}(\tau)\big|d\tau \\
+\int_s^t\int_\Omega e^{\mu\tau}\big|\partial_tf(u^{\eps\delta}(\tau))\big|d\tau+\int_s^t
\int_\Omega e^{\mu\tau}\partial_u f(u^{\eps\delta}(\tau))\, \big|\partial_t u^{\eps\delta}(\tau)\big|\,d\tau.
\end{multline}

By the assumption (\ref{Phif}), $\partial_u f\le-\mu<0$ and setting $R$ such that,
\begin{equation}\label{38}
|u^{\eps\delta}(x,t)|\le R\qquad\forall t>0,\ x\in\Omega,
\end{equation}
uniformly in $\eps$ and $\delta$, we deduce from (\ref{37}) that
\begin{equation}\label{39}
e^{\mu t}\nu(t)\le e^{\mu s}\nu(s)+\int_s^t e^{\mu\tau}\xi_R(\tau)\,d\tau,\qquad t>s\ge 0,
\end{equation}
where we set $\nu(t)=\displaystyle\int_\Omega|\partial_t u^{\eps\delta}(t)|$ and $\xi_R(t)=\displaystyle\int_\Omega\sup_{|u|\le R}|\partial_t f(u)|$.

>From (\ref{39}) and with an elementary estimate of its last term, we obtain the following inequality for any $t>0$, $\sigma\ge 0$,
\begin{equation*}
\nu(t+\sigma)\le e^{-\mu t}\nu(\sigma)+C_\mu\sup_{\sigma<s<t+\sigma}\int_s^{s+1}\xi_R(\tau)\,d\tau
\end{equation*}
with $C_\mu=(1-e^{-\mu})^{-1}+1$. Consequently, from (\ref{xiR}), we conclude
\begin{equation*}
\nu(t+\sigma)\le e^{-\mu t}\big(\nu(0)+C_\mu\, C_R\big)+ C_\mu\sup_{\sigma<s<t+\sigma}\int_s^{s+1}\xi_R(\tau)\,d\tau,
\end{equation*}
and the right hand side vanishes, as $t\rightarrow\infty$, uniformly in $\eps$ and $\delta$, provided (\ref{38}) holds. But the $L^\infty$ estimate
of $u$ globally in time holds with $R=R(\Lambda,diam(\Omega))$ in the case of assumption (\ref{gu})  or with $R=\max\{M,\|u_0\|_{L^\infty(\Omega)}\}$,
 where $M$  is given by the assumption
(\ref{supsubsol}) in the second case. Notice that, when (\ref{supsubsol}) is satisfied, by the weak maximum principle (see, for instance Theorem 1.6 of
\cite{ChipotRodrigues1988}), $R$ is a supersolution and $-R$ a subsolution of (\ref{aprox}).

Therefore, it is clear that, by compactness, we may take a subsequence $t_n\rightarrow\infty$ and a
function $u^\infty\in W^{1,\infty}_0(\Omega)$ such that ($0<\alpha<1$)
\begin{equation}\label{48}
u(t_n)\underset{n}{\longrightarrow}u^\infty\qquad\text{ in }W^{1,\infty}_0(\Omega)\text{ weak\,-\,}*\ \text{ and in }C^{0,\alpha}(\bar\Omega) \text{
strong},
\end{equation}
\begin{equation}\label{49}
\partial_tu(t_n)\underset{n}{\longrightarrow}0\qquad\text{ in }M(\Omega).
\end{equation}

For an arbitrary measurable subset $\omega\subset\Omega$, we have

\begin{equation*}
\int_\omega\big|\nabla u^\infty\big|\le\liminf_n\int_\omega\big|\nabla u(t_n)\big|\le\lim_n\int_\omega G(u(t_n))=\int_\omega G(u^\infty),
\end{equation*}
we obtain that
\begin{equation*}
\big|\nabla u^\infty(x)\big|\le G(x,u^\infty(x))\qquad\mbox{ a.e. in }x\in\Omega,
\end{equation*}
and so $u^\infty\in\K_{G(u^\infty)}$.

Consider the quasivariational inequality (\ref{iqv}) at $t=t_n$
\begin{multline}
\label{50}
\int_\Omega\partial_tu(t_n)\,\big(v-u(t_n)\big)+\int_\Omega\boldsymbol{\Phi}\big(u(t_n)\big)\cdot\nabla\big(v-u(t_n)\big)\\
\ge\int_\Omega f(t_n,u(t_n))\,\big(v-u(t_n)\big),\qquad\forall v\in\K_{G(u(t_n))}.\end{multline}

Given $\sigma\in(0,1)$ and $w^\infty$ an arbitrary function of $\K_{G(u^\infty)}$ in $L^\infty(\Omega)$, as
$G(u(t_n))\underset{n}{\longrightarrow}G(u^\infty)$ uniformly, there exists $n_0\in\N$ such that, for $n\ge n_0$ we have a.e. ${x\in\Omega}$ 
\begin{equation*}
 G(u^\infty)\le\frac1{(1-\sigma)}G(u(t_n)).
 \end{equation*}

Defining
\begin{equation*}
v^\sigma(x)=(1-\sigma)\,w^\infty(x),
\end{equation*}
we have
\begin{equation*}
\big|\nabla v^\sigma\big|=\,(1-\sigma)\,\big|\nabla w^\infty\big|\le\,(1-\sigma)\,G(u^\infty)\le\, G(u(t_n)).
\end{equation*}
which means that $v^\sigma\in\K_{G(u(t_n))}$, for all $n\ge n_0$. So, taking $v=v^\sigma$ in (\ref{50}) and letting $t_n\rightarrow\infty$  we obtain
\begin{equation*}
\int_\Omega\boldsymbol{\Phi}(u^\infty)\cdot\nabla\big((1-\sigma) w^\infty-u^\infty\big)\ge\int_\Omega f^\infty
(u^\infty)\big((1-\sigma)w^\infty-u^\infty\big),
\end{equation*}
since by the assumption (\ref{fff}) and the convergence (\ref{48}) we have
\begin{equation*}
f(t_n,u(t_n))\underset{n}{\longrightarrow}f^\infty(u^\infty)\qquad\text{ in }L^p(\Omega),\qquad 1\le p<\infty.
\end{equation*}

Finally, letting $\sigma\rightarrow 0$,  $u^\infty$ solves problem (\ref{iqvs}), as $w^\infty$ is an arbitrary function of $\K_{G(u^\infty)}$.
\hfill{$\square$}

\vspace{3mm}

\noindent{\em Proof of Proposition \ref{asymp}.}

\noindent Let $\vartheta$ be the solution of problem (\ref{iv}) and $\vartheta^\infty$ the solution of problem
(\ref{ivinfty}).

Recall that  $v=\vartheta+\zeta^2s_\zeta(\vartheta^\infty-\vartheta)$, $0<\zeta<1$, belongs to $\K_G$. So, using $v$ as test function in (\ref{iv}) and calling $\overline\vartheta=\vartheta^\infty-\vartheta$, we have
\begin{equation*}
\zeta^2\int_\Omega\partial_t\vartheta(t)\,s_\zeta(\overline\vartheta)+\zeta^2\int_\Omega s_\zeta'
(\overline\vartheta)\boldsymbol\Phi(\vartheta)\cdot\nabla\overline\vartheta
\ge \zeta^2\int_\Omega f(\vartheta)\,s_\zeta(\overline\vartheta).
\end{equation*}

Using $w=\vartheta^\infty+\zeta^2s_\zeta(\overline\vartheta)\in\K_G$ as test function in (\ref{ivinfty}), we get
\begin{equation*}
-\zeta^2\int_\Omega s_\zeta'(\overline\vartheta)\boldsymbol\Phi(\vartheta^\infty)\cdot\nabla\overline\vartheta
\ge -\zeta^2\int_\Omega f^\infty(\vartheta^\infty)\,s_\zeta(\overline\vartheta).
\end{equation*}

Summing  the above inequalities we get
\begin{equation*}
\int_\Omega\partial_t\overline\vartheta\,s_\zeta(\overline\vartheta)+
\int_\Omega s_\zeta'(\overline\vartheta)\big(\boldsymbol\Phi(\vartheta^\infty)-\boldsymbol\Phi(\vartheta)\big)
\cdot\nabla\overline\vartheta
\le \int_\Omega\big( f^\infty(\vartheta^\infty)-f(\vartheta)\big)\,s_\zeta(\overline\vartheta).
\end{equation*}

After multiplication by $e^{\mu\tau}$ and integration in order to $\tau$, between $\sigma$ and $t$, since
\begin{equation*}
\int_\sigma^te^{\mu\tau}\frac{d}{d\tau}\int_\Omega\partial_t\overline\vartheta\,s_\zeta(\overline\vartheta)\,d\tau=
e^{\mu t}\int_\Omega S_\zeta(\overline\vartheta(t))-e^{\mu \sigma}\int_\Omega S_\zeta(\overline\vartheta(\sigma))-
\mu\int_\sigma^t\int_\Omega e^{\mu\tau} S_\zeta(\overline\vartheta)\,d\tau,
\end{equation*}
we have
\begin{multline*}
\int_\Omega S_\zeta(\overline\vartheta(t))\le e^{\mu(\sigma-t)}\int_\Omega S_\zeta(\overline\vartheta(\sigma))+\int_\sigma^t\int_\Omega
e^{\mu\tau} \big(\mu S_\zeta(\overline\vartheta)-(f^\infty(\vartheta)-f^\infty(\vartheta^\infty))s_\zeta(\overline\vartheta)\big)\,d\tau\\
+\int_\sigma^te^{\mu(\tau-t)}\int_\Omega s_\zeta'(\overline\vartheta)\big(\boldsymbol\Phi(\vartheta^\infty)-
\boldsymbol\Phi(\vartheta)\big)\cdot\nabla\overline\vartheta\,d\tau
+\int_\sigma^t e^{\mu(\tau-t)}\int_\Omega\big| f(\vartheta^\infty)-f(\vartheta)\big|\,d\tau.
\end{multline*}

Using the strict decreasing property (14), we find
\begin{align*}
\int_\sigma^t\int_\Omega
e^{\mu\tau} \big(\mu S_\zeta(\overline\vartheta)&-(f^\infty(\vartheta)-f^\infty(\vartheta^\infty))s_\zeta(\overline\vartheta)\big)\,d\tau=
\int_\sigma^t\int_\Omega
e^{\mu\tau} \big(\mu S_\zeta(\overline\vartheta)+\partial_u f^\infty(\xi)\overline\vartheta s_\zeta(\overline\vartheta))\,d\tau\\
&\le\int_\sigma^t\int_\Omega e^{\mu\tau} \mu\big(S_\zeta(\overline\vartheta)-\overline\vartheta s_\zeta(\overline\vartheta)\big)\,d\tau
\underset{\zeta\rightarrow 0}{\longrightarrow} 0 .
\end{align*}
On the other hand,
\begin{equation*}
\Big|\int_\sigma^te^{\tau-t}\int_\Omega s_\zeta'(\overline\vartheta)\big(\boldsymbol\Phi(\vartheta^\infty)-
\boldsymbol\Phi(\vartheta)\big)\Big|\,d\tau\le\int_\Omega L_{\boldsymbol\Phi}\big|\overline\vartheta\big|\,s_\zeta'(\overline\vartheta)\,
\big|\nabla\overline\vartheta\big|\,d\tau,
\end{equation*}
where $L_{\boldsymbol\Phi}$ is a Lipschitz constant, also converges to zero when $\zeta\rightarrow 0$, since
$\displaystyle\lim_{\zeta\rightarrow 0}\tau\,s'_\zeta(\tau)=0$.
Therefore, we obtain
\begin{equation*}
\int_\Omega \big|\overline\vartheta(t+\sigma)\big|\le e^{-\mu t} \int_\Omega
\big|\overline\vartheta(\sigma)\big|+\int_\sigma^{t+\sigma} e^{\mu(\tau-t-\sigma)}\int_\Omega\big| f(\vartheta^\infty)-f(\vartheta)\big|\,d\tau.
\end{equation*}

Arguing as in \cite{ChipotRodrigues1988}, page 283, we have
\begin{equation*}
\big\|\vartheta(t)-\vartheta^\infty\big\|_{L^1(\Omega)}\le e^{-\mu t}\big\|\vartheta(\sigma)-\vartheta^\infty\big\|_{L^1(\Omega)}
+C\sup_{\sigma<\tau<t}\eta_{_M}(\tau)|,
\end{equation*}
and so $\vartheta(t)\underset{t\rightarrow\infty}{\longrightarrow}\vartheta^\infty$ in $L^1(\Omega)$, since $\eta_{_M}(t)\underset{t
\rightarrow\infty}{\longrightarrow}0$.
As $\|\nabla(\vartheta(t)-\vartheta^\infty)\|_{L^\infty(\Omega)}\le 2\Lambda$  then
$\vartheta(t)\underset{t\rightarrow\infty}{\longrightarrow}\vartheta^\infty$ in $C^{0,\alpha}(\bar\Omega)$,  for $0\le\alpha<1$, by compactness.

If, in addition, $\eta_{_M}(t)=O(e^{-\gamma t})$, $\gamma>0$, when $t\rightarrow\infty$,
then $\|\vartheta-\vartheta^\infty\|_{L^1(\Omega)}=O(e^{-\nu t})$, where $\nu=\min\{\mu,\gamma\}$ and,
by a  Gagliardo-Niremberg interpolation inequality, we conclude as in \cite{Rodrigues1987}, page 200, for any $0\le\alpha<1$,
\begin{equation*}
\|\vartheta(t)-\vartheta^\infty\|_{C^{0,\alpha}(\bar\Omega)}=O\big(e^{-\big(\frac{1-\alpha}{n+1}\big)\nu t}\big),\qquad\text{as }t\rightarrow\infty.
\end{equation*}
\hfill{$\square$}

\vspace{5mm}

%

\def\ocirc#1{\ifmmode\setbox0=\hbox{$#1$}\dimen0=\ht0
    \advance\dimen0 by1pt\rlap{\hbox to\wd0{\hss\raise\dimen0
    \hbox{\hskip.2em$\scriptscriptstyle\circ$}\hss}}#1\else
    {\accent"17 #1}\fi}

\begin{small}

\begin{tabular}{ll}
\noindent Jos\'e Francisco Rodrigues&\hspace{0,5cm}Lisa Santos\\
CMAF and Department of Mathematics,&\hspace{0,5cm}CMAT / Dept. of Mathematics and Applications\\
University of Lisbon,&\hspace{0,5cm}University of Minho,\\
Av. Prof. Gama Pinto, 2,&\hspace{0,5cm}Campus de Gualtar,\\
1649-003 Lisboa, Portugal&\hspace{0,5cm}4710-057 Braga, Portugal\\
rodrigue@fc.ul.pt&\hspace{0,5cm}lisa@math.uminho.pt
\end{tabular}

\end{small}

\end{document}